\documentclass{amsproc}

\usepackage[english]{babel}
\usepackage[cp1250]{inputenc}
\usepackage{amsmath,amssymb,amsfonts}
\usepackage{url}
\usepackage[dvipsnames,usenames]{color}
\usepackage{array}
\usepackage[pdftex]{graphicx}
\usepackage{tikz} 

\newcolumntype{L}{>{\displaystyle}l}
\newcolumntype{C}{>{\displaystyle}c}
\newcolumntype{R}{>{\displaystyle}r}

\renewcommand{\tfrac}{\genfrac{}{}{}1}

\newcommand{\R}{\mathbb R}
\newcommand{\N}{\mathbb N}
\newcommand{\Z}{\mathbb Z}
\newcommand{\C}{\mathbb C}

\newcommand{\cont}{\mathcal C}

\newcommand{\D}{\mathbb D}
\renewcommand{\Im}{\mathrm{Im}}
\renewcommand{\Re}{\mathrm{Re}}
\def\E{{\mathrm{e}}}

\def\til{\widetilde}

\newcommand{\pois}{\mathcal P}

\def\I{\mathfrak{i}}

\newcommand{\diff}{\mathrm{d}}
\renewcommand{\bar}{\overline}

\newcommand{\och}{\quad\mathrm{and}\quad}

\newcommand{\myref}[1]{$(\ref{#1})$}

\theoremstyle{definition}
\newtheorem{definicija1}{Definition}[section]
\newtheorem{primer1}[definicija1]{Example}
\newtheorem{opomba1}[definicija1]{Remark}

\theoremstyle{plain}

\newtheorem{izrek1}[definicija1]{Theorem}
\newtheorem{trditev1}[definicija1]{Proposition}
\newtheorem{posledica1}[definicija1]{Corollary}

\numberwithin{equation}{section}

\begin{document}

\title[Herglotz-Nevanlinna functions in two variables]{A characterization of Herglotz-Nevanlinna functions in two variables via integral representations}

\author{Annemarie Luger}
\address{Annemarie Luger, Department of Mathematics, Stockholm University, 106 91 Stockholm, Sweden}
\curraddr{}
\email{luger@math.su.se}

\author{Mitja Nedic}
\address{Mitja Nedic, Department of Mathematics, Stockholm University, 106 91 Stockholm, Sweden}
\curraddr{}
\email{mitja@math.su.se}
\thanks{The authors are supported by the SSF applied mathematics grant nr. AM13-0011.}


\date{24-3-2016}

\begin{abstract}
We derive an integral representation for Herglotz-Nevanlinna functions in two variables which provides a complete characterization of this class  in terms of a real number, two non-negative numbers and a positive measure satisfying certain conditions. Further properties of the representing measures are discussed.
\end{abstract}

\maketitle

\section{Introduction}

Herglotz-Nevanlinna functions in one variable are functions analytic in the upper half plane having non-negative imaginary part. This class has been very well studied during the last century and has proven very useful in many applications. It seems natural to consider corresponding functions in several variables, i.e. analytic functions that have non-negative imaginary part if all variables lie in the upper half plane. From applications point of view such functions are very interesting when considering linear passive systems with several parameters. 

In the treatment of Herglotz-Nevanlinna functions in one variable, a very strong tool is the classical result that these functions can be characterized via integral representations, cf. Theorem \ref{thm21}. It is hence a natural question to ask for a corresponding representation for functions in several variables. 

Already almost 50 yeast ago the form of such integral representations was suggested by Vladimirov, see e.g. \cite{vladimirov1}.  Using a very heavy machinery of classical distribution theory, it is shown that every Herglotz-Nevanlinna function can be written in such a form, but it is a priory assumed that the measure appearing in the formula is the boundary measure of the function. Thus the drawback in these results is that they do not specify the properties of the measure and hence cannot provide a characterization as in the case of only one variable. In \cite{vladimirov2} the authors use a different approach in order to find a characterization, however, the obtained representation becomes much more involved. In view of the present result one can in fact say that both the representation and the conditions on the measure are \emph{too} complicated since it turns out that many terms there actually vanish or simplify radically, see Remark \ref{vlad} for more details. 

Recently the question of  a characterization was taken up again by Agler, McCarthy and Young in \cite{agler2}, see also Agler, Tully-Doyle and Young in \cite{agler}. They found a characterization via operator representations, however, only for certain subclasses of functions satisfying an asymptotic condition. 
 
In the present paper, we solve the characterization problem for the whole class of Herglotz-Nevanlinna functions in two variables by deriving an integral representation together with conditions on the representing measure. Even if these representations are of the same form as in \cite{vladimirov1}, our result contains considerably more information, since we obtain a full (but simple) description of all representing measures. Moreover, the proof is shorter and  more elementary in the sense that is uses only Theorems \ref{thm21} and \ref{04thm2} which are in their essence built upon Cauchy's integral formula and Helly's selection principle. We use - as in the classical proof for the one dimensional case - a corresponding result for the polydisk, and then a transformation of variables. However, unlike in one variable the terms arising from the boundary of the area of integration are quite delicate and need very careful treatment in order to simplify the representation to the desired form.

It appears that the requirements on the representing measures have quite strong consequences, which are discussed in Section \ref{measure4}. 

\section{Notations and a brief recap of known results}\label{2}

As usual $\D$ denotes the unit disk in the complex plane while $\C^+$ denotes the upper and  $\C_+$ denotes the right half-plane. Throughout this paper we will use the convention that $z$ denotes the complex variable that lies in the upper half-plane while $w$ denote the variable that lies in the disk. We recall also the fact that the unit disk and the upper half-plane are biholomorphic. One map achieving this is $\varphi\colon \C^+ \to \D$ defined as $\varphi\colon z \mapsto \tfrac{z-\I}{z+\I}$. Its inverse is then given as $\varphi^{-1}\colon w \mapsto\I\tfrac{1+w}{1-w}$. Note also that $\varphi$ is a bijection between the sets $\R$ and $S^1 \setminus \{1\}$.

It is often convenient to consider Herglotz-Nevanlinna functions that do not attain real values. In particular, every such function $q\colon \C^+ \to \C^+$ then uniquely determines a function $f\colon \D \to \C_+$ with respect to the biholomorphisms $\varphi$ and $\cdot \I$, as elaborated by the diagram in Figure \ref{fig1}. The converse also holds; a function $f$ uniquely determines $q$ with respect to the same biholomorphisms. It can be shown that the only Herglotz-Nevanlinna functions that are excluded form this correspondence are in fact real-constant functions.

\begin{figure}[!h]\centering

\begin{tikzpicture}[node distance=2cm, auto]
\node (up1) {$\C^+$};
\node (up2) [right of=up1] {$\C^+$};
\node (disk) [below of=up1] {$\D$};
\node (right) [below of=up2] {$\C_+$};
\draw[->] (up1) to node {$q$} (up2);
\draw[->] ([xshift=-2pt]up1.south) to node [swap] {$\varphi$} ([xshift=-2pt]disk.north);
\draw[->] ([xshift=2pt]disk.north) to node [swap] {$\varphi^{-1}$} ([xshift=2pt]up1.south);
\draw[->] (disk) to node [swap] {$f$} (right);
\draw[->] ([xshift=-2pt]right.north) to node {$\cdot \I$} ([xshift=-2pt]up2.south);
\draw[->] ([xshift=2pt]up2.south) to node {$\cdot \I^{-1}$} ([xshift=2pt]right.north);
\end{tikzpicture}

\caption{The relationship between $q$ and $f$.}
\label{fig1}
\end{figure}

We recall now the integral representation theorem due to Nevanlinna \cite{nevan1}, which was presented in its current form by Cauer \cite{cauer}.

\begin{izrek1}[Nevanlinna]\label{thm21}
A function $q\colon \C^+ \to \C$ is a Herglotz-Nevanlinna function if and only if $q$ can be written as
\begin{equation}\label{eqthm21}
q(z) = a + b z + \frac{1}{\pi}\int_{-\infty}^\infty\left(\frac{1}{t-z} - \frac{t}{1 + t^2}\right)\diff\mu(t)
\end{equation}
where $a \in \R$, $b \geq 0$ and $\mu$ is a positive Borel measure on $\R$ satisfying
\begin{equation}\label{eqthm22}
\int_{-\infty}^\infty\frac{1}{1 + t^2}\diff\mu(t) < \infty.
\end{equation}
\end{izrek1}

\begin{opomba1}
Moreover, $a$, $b$ and $\mu$ are unique with these properties.\end{opomba1}

The importance and beauty of the theorem is that it gives a complete characterization of $q$ in terms of the numbers $a$ and $b$ and the measure $\mu$. But it also provides a tool for handling Herglotz-Nevanlinna functions. We mention the following property that will be of use to us further on.

\begin{trditev1}\label{prop22}
Let $q$ be a Herglotz-Nevanlinna function. Then the non-tangential limit
\begin{equation}\label{eqprop21}
\lim\limits_{z \hat{\to} \infty}\frac{q(z)}{z} = b,
\end{equation}
where $b\geq0$ is the number that appears in representation \myref{eqthm21}.
\end{trditev1}

Recall that $z \hat{\to} \infty$ is a shorthand notation for $|z| \to \infty$ in the Stoltz domain $\{z \in \C^+~|~\theta \leq \arg(z) \leq \pi-\theta\}$ for any $\theta \in (0,\tfrac{\pi}{2}]$.

Let us denote by $\C^{+2} := \{z = (z_1,z_2) \in \C^2~|~\Im[z_1] > 0, \Im[z_2] > 0\}$ the poly-upper half-plane in $\C^2$. Our main object of interest is the following class of functions in two variables.

\begin{definicija1}\label{def1}
A holomorphic function $q\colon \C^{+2} \to \C$ with non-negative imaginary part is called a \emph{Herglotz-Nevanlinna function} (\emph{in two variables}).
\end{definicija1}

In the situation of one variable, a standard proof of Theorem \ref{thm21} uses the Riesz-Herglotz theorem, see e.g. \cite{cauer}, that gives an integral representation for functions on the unit disk with positive real part. It is then possible to use the biholomorphisms discussed earlier to return to functions defined on the upper half-plane.

In order to apply the same strategy in several variables we use a generalization of the Riesz-Herglotz theorem by the Kor\'{a}nyi-Puk\'{a}nszky, \cite{koranyi}, that completely characterizes functions defined on the unit polydisk in $\C^n$ that have positive real part. It seems that Vladimirov has  independently the same result in \cite{vladimirov3}, which is used in \cite{vladimirov2}. Here we present the theorem only for $n=2$ and with slightly different notation that is more inclined towards our purpose of giving a representation of Herglotz-Nevanlinna functions. 

\begin{izrek1}\label{04thm2}
A function $f$ on the unit polydisk $\D^2$ is holomorphic and has non-negative real part if and only if $f$ can be written as
\begin{multline}\label{04eq4}
f(w_1,w_2) = \I~\Im[f(0,0)] \\
+ \frac{1}{4\pi^2}\iint_{[0,2\pi)^2}\left(\frac{2}{(1-w_1\E^{-\I s_1})(1-w_2\E^{-\I s_2})}-1\right)\diff\nu(s_1,s_2)
\end{multline}
where $\nu$ is a finite positive Borel measure on $[0,2\pi)^2$ satisfying the condition that
\begin{equation}\label{04eq5}
\iint_{[0,2\pi)^2}\E^{\I m_1 s_1}\E^{\I m_2 s_2}\diff\nu(s_1,s_2) = 0
\end{equation}
for every pair of indices $m_1,m_2 \in \Z$ satisfying $m_1m_2 < 0$.
\end{izrek1}

\section{The theorem in two variables}

Before presenting the main theorem we introduce some notation. Denote by $K$ the kernel function depending on $(z_1,z_2) \in \C^{+2}$ and $(t_1,t_2) \in \R^2$, defined as
$$K\big((z_1,z_2),(t_1,t_2)\big): = -\frac{\I}{2}\:\!\!\left(\frac{1}{t_1-z_1}-\frac{1}{t_1+\I}\right)\:\left(\frac{1}{t_2-z_2} - \frac{1}{t_2+\I}\right) + \frac{\I}{(1+t_1^2)(1+t_2^2)}.$$
We will also need the Poisson kernel of $\C^{+2}$, which we prefer to write using complex coordinates as
$$\pois\big((z_1,z_2),(t_1,t_2)\big) := \frac{\Im[z_1]}{|t_1 - z_1|^2}\:\frac{\Im[z_2]}{|t_2 - z_2|^2}.$$
Note that $\pois > 0$ for any $(z_1,z_2) \in \C^{+2}$ and any $(t_1,t_2) \in \R^2$.

The main result of this paper is as follows.

\begin{izrek1}\label{thm5}
A function $q\colon \C^{+2} \to \C$ is a Herglotz-Nevanlinna function in two variables  if and only if $q$ can be written as
\begin{equation}\label{eqthm51}
q(z_1,z_2) = a + b_1 z_1 + b_2 z_2 + \frac{1}{\pi^2}\iint_{\R^2}K\big((z_1,z_2),(t_1,t_2)\big)\diff\mu(t_1,t_2)
\end{equation}
where $a \in \R$, $b_1,b_2 \geq 0$, and $\mu$ is a positive Borel measure on $\R^2$ satisfying the growth condition
\begin{equation}\label{eqthm52}
\iint_{\R^2}\frac{1}{(1 + t_1^2)(1+t_2^2)}\diff\mu(t_1,t_2) < \infty
\end{equation}
and the Nevanlinna condition
\begin{equation}\label{eqthm53}
\iint_{\R^2}\Re\left[\left(\frac{1}{t_1-\bar{z_1}}-\frac{1}{t_1+\I}\right)\left(\frac{1}{t_2-z_2}-\frac{1}{t_2-\I}\right)\right]\diff\mu(t_1,t_2) \equiv 0
\end{equation}
for all $(z_1,z_2) \in \C^{+2}$.
\end{izrek1}

\begin{opomba1}
Moreover, the numbers $a$, $b_1$, $b_2$ and the measure $\mu$ are uniquely determined as it is shown in Corollary \ref{posfour} and Proposition \ref{pos4}, respectively.
\end{opomba1}

\begin{opomba1}
Observe that  for functions of one variable there is no analogue  to the Nevanlinna condition \eqref{eqthm53}.
\end{opomba1}

\proof
Let us assume first that $q$ is a Herglotz-Nevanlinna function. We first consider the possibility that $q$ attains a real value. Then there exists a point $(\zeta_1,\zeta_2) \in \C^{+2}$ such that $\Im[q(\zeta_1,\zeta_2)] = 0$. Since $q$ is a Herglotz-Nevanlinna function it is holomorphic and its imaginary part $\Im[q]\geq0$ is therefore pluriharmonic. It follows now from the maximum principle for pluriharmonic functions that $\Im[q] \equiv 0$ on $\C^{+2}$, and hence the function $q$ admits  a representation of the form \myref{eqthm51} with $a = q(\zeta_1,\zeta_2)$, $b_1 = b_2 = 0$ and $\mu \equiv 0$. Thus the theorem holds in this case.

We may now restrict ourselves to the case when $q$ does not attain a real value. Then there exists a function $f$ on $\D^2$ with positive real part such that
$$q(z_1,z_2) = \I f(\varphi(z_1),\varphi(z_2)),$$
where $\varphi$ is given as in Section \ref{2}.
Using representation \myref{04eq4} of the function $f$ yields
\begin{eqnarray*}
q(z_1,z_2) & = &  -\Im[f(0,0)] \\
~& ~& + \frac{\I}{4\pi^2}\iint_{[0,2\pi)^2}\left(\frac{2}{(1-\varphi(z_1)\E^{-\I s_1})(1-\varphi(z_2)\E^{-\I s_2})}-1\right)\diff\nu(s_1,s_2).
\end{eqnarray*}
We obtain the first term in  representation \myref{eqthm51} by setting $a := -\Im[f(0,0)] \in \R$.

\begin{figure}[!ht]
\centering
\includegraphics[height=55mm]{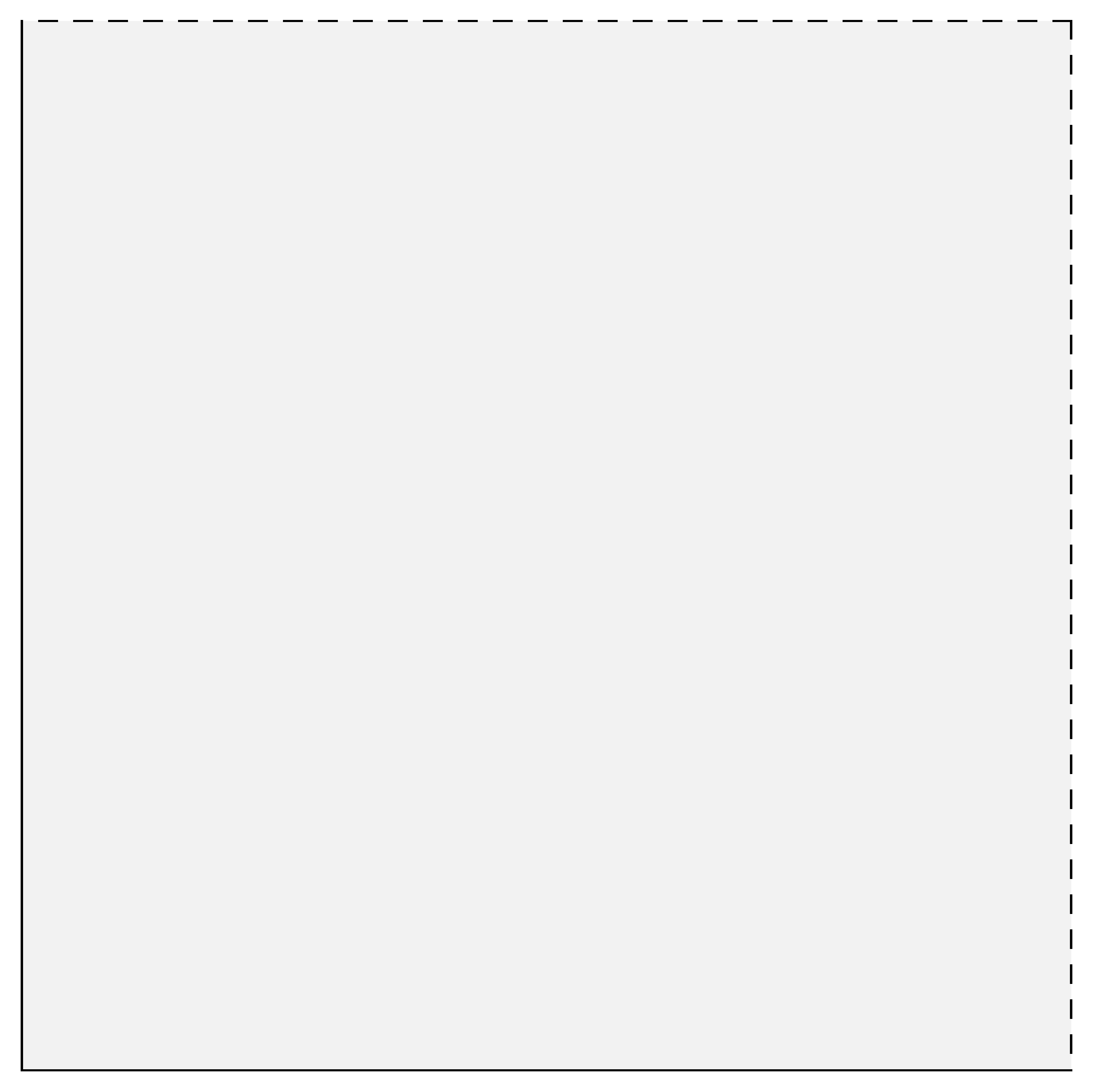}
\caption{The area of integration.}
\label{fig2}
\end{figure}

Before transforming the area of integration to $\R^2$ we divide the integral over $[0,2\pi)^2$, which is shown in Figure \ref{fig2}, into four parts and investigate each part separately.

Considering first the integral over the open square $(0,2\pi)^2$ we do a change of variables where $\E^{\I s_i} = \tfrac{t_i - \I}{t_i + \I}$ for $i=1,2$. The area of integration thus transforms into $\R^2$ and the measure $\nu$ transforms into a measure $\mu$ related by the chosen change of variables as
$$\diff\nu(s_1,s_2) = \frac{4}{(1+t_1^2)(1+t_2^2)}\diff\mu(t_1,t_2).$$
As an immediate consequence of this transformation we see that the measure $\mu$ satisfies condition \myref{eqthm52} since $\nu$ is a finite measure and
$$\iint_{\R^2}\frac{1}{(1 + t_1^2)(1+t_2^2)}\diff\mu(t_1,t_2) = \frac{1}{4}\iint_{(0,2\pi)^2}\diff\nu(s_1,s_2) < \infty.$$
The integral thus becomes
\begin{multline*}
\frac{\I}{4\pi^2}\iint_{(0,2\pi)^2}\left(\frac{2}{(1-\varphi(z_1)\E^{-\I s_1})(1-\varphi(z_2)\E^{-\I s_2})}-1\right)\diff\nu(s_1,s_2) \\
= \frac{1}{\pi^2}\iint_{\R^2}K\big((z_1,z_2),(t_1,t_2)\big)\diff\mu(t_1,t_2)
\end{multline*}
where the equality between the two expressions comes exclusively from symbolic manipulations of the first term along with the discussed change of variables. This gives us the integral term of representation \myref{eqthm51}.

We now consider the part of the integral that runs over one side of the square, namely   $\{0\} \times (0,2\pi)$. Let us denote $\diff\nu_1(s_2) := \diff\nu(0,s_2)$ and let $\mu_1$ be a measure on $\R$ related to $\nu_1$ as
$$\diff\nu_1(s_2) = \frac{2}{1+t_2^2}\diff\mu_1(t_2).$$
Thus
\begin{multline*}
\frac{\I}{4\pi^2}\int_{\{0\} \times (0,2\pi)}\left(\frac{2}{(1-\varphi(z_1)\E^{-\I s_1})(1-\varphi(z_2)\E^{-\I s_2})}-1\right)\diff\nu(s_1,s_2) \\
=  \frac{1}{2\pi^2}\int_{-\infty}^{\infty}\left(\frac{z_1+\I}{2\I}\:\frac{1+t_2z_2}{t_2-z_2} + \frac{z_1-\I}{2\I}\:\frac{1+t_2 \I}{t_2-\I}\right)\frac{1}{1+t_2^2}\diff\mu_1(t_2).
\end{multline*}
For simplicity we introduce  the function $q_1$  defined as
$$q_1(\zeta) := \frac{1}{\pi}\int_{-\infty}^{\infty}\frac{1+t \zeta}{t-\zeta}\:\frac{1}{1+t^2}\diff\mu_1(t),$$
then the above integral can be written as 
$$\frac{1}{2\pi}\left(\frac{z_1 + \I}{2\I}q_1(z_2) + \frac{z_1 - \I}{2\I}q_1(\I)\right).$$
We note that $q_1$ is a Herglotz-Nevanlinna function in one variable with the numbers $a$ and $b$ from representation \myref{eqthm21} both equal to 0 and the measure equal to $\mu_1$.

An analogous procedure for the other side of the square gives that
\begin{multline*}
\frac{\I}{4\pi^2}\int_{(0,2\pi) \times \{0\}}\left(\frac{2}{(1-\varphi(z_1)\E^{-\I s_1})(1-\varphi(z_1)\E^{-\I s_2})}-1\right)\diff\nu(s_1,s_2) \\
= \frac{1}{2\pi}\left(\frac{z_2 + \I}{2\I}q_2(z_1) + \frac{z_2 - \I}{2\I}q_2(\I)\right)
\end{multline*}
where the  function $q_2$ are defined in an analogous way as in the previous case.

Finally, integration over the corner point  $\{0\} \times \{0\}$ gives that
\begin{multline*}\frac{\I}{4\pi^2}\int_{\{0\} \times \{0\}}\left(\frac{2}{(1-\varphi(z_1)\E^{-\I s_1})(1-\varphi(z_2)\E^{-\I s_2})}-1\right)\diff\nu(s_1,s_2) \\
=  \frac{\I}{4\pi^2}\left(\frac{2}{(2\I)^2}(z_1+\I)(z_2+\I) - 1\right)\nu(\{(0,0)\}).
\end{multline*}

Hence we have so far arrived at a representation of the function $q$ of the form
\begin{eqnarray}\label{eq52}
q(z_1,z_2) & = & a + \frac{\I}{4\pi^2}\left(\frac{2}{(2\I)^2}(z_1+\I)(z_2+\I) - 1\right)\nu(\{(0,0)\}) \\
~& ~& +\frac{1}{2\pi}\left(\frac{z_2 + \I}{2\I}q_2(z_1) + \frac{z_2 - \I}{2\I}q_2(\I)\right) \nonumber \\
~& ~ & + \frac{1}{2\pi}\left(\frac{z_1 + \I}{2\I}q_1(z_2) + \frac{z_1 - \I}{2\I}q_1(\I)\right) \nonumber \\ 
 ~ & ~ & + \frac{1}{\pi^2}\iint_{\R^2}K\big((z_1,z_2),(t_1,t_2)\big)\diff\mu(t_1,t_2). \nonumber
\end{eqnarray}
While the first and last part are already as desired, we still have to show that the middle three terms indeed give the two linear terms from representation \myref{eqthm51}. This will be done by showing that the functions $q_1$ and $q_2$ actually are of a very particular form. 

Let  $\alpha \in \C^+$ and consider the function $\til{q}_{1}$ defined by fixing the second variable  $\til{q}_{1}(z) := q(z,\alpha)$ for $z \in \C^+$, which is a Herglotz-Nevanlinna function in one variable. In view of Proposition \ref{prop22} we consider the following non-tangential limit 
\begin{eqnarray}
\lim\limits_{z \hat{\to} \infty}\frac{\til{q}_1(z)}{z} & = &  \lim\limits_{z \hat{\to} \infty}\frac{a}{z} + \lim\limits_{z \hat{\to} \infty}\frac{\I}{(2\pi)^2}\frac{1}{z}\left(\frac{2}{(2\I)^2}(z+\I)(\alpha+\I) - 1\right)\nu(\{(0,0)\}) \nonumber\\
~& ~ & +  \lim\limits_{z \hat{\to} \infty}\frac{1}{2\pi}\left( \frac{\alpha + \I}{2\I z}q_2(z) + \frac{\alpha - \I}{2\I z}q_2(\I)\right) \nonumber\\
~& ~ & +  \lim\limits_{z \hat{\to} \infty}\frac{1}{2\pi}\left( \frac{z + \I}{2\I z}q_1(\alpha) + \frac{z - \I}{2\I z}q_1(\I)\right)  \nonumber\\
~& ~ & + \lim\limits_{z \hat{\to} \infty}\frac{1}{\pi^2}\iint_{\R^2}\frac{K\big((z,\alpha),(t_1,t_2)\big)}{z}\diff\mu(t_1,t_2)   \nonumber\\
~ & = & - \frac{\I}{4\pi^2}\frac{1}{2}(\alpha+\I)\nu(\{(0,0)\}) +  \frac{1}{2\pi}\:\frac{q_1(\alpha)+q_1(\I)}{2\I}.\label{lim}
\end{eqnarray}
Here we used  that
$$\lim\limits_{z \hat{\to} \infty}\frac{q_2(z)}{z} = 0$$
by Proposition \ref{prop22} and that the last term vanishes as the interchange of the limit and integral  is valid as the assumptions of Lebesgue's dominated convergence theorem are satisfied.

Again Proposition \ref{prop22} now implies that the above limit in \eqref{lim} is non-negative, i.e. for $\alpha \in \C^+$ it holds 
\begin{equation}\label{eq563}
\frac{1}{2\pi}\:\frac{q_1(\alpha)+q_1(\I)}{2\I} - \frac{\I}{4\pi^2}\frac{1}{2}(\alpha+\I)\nu(\{(0,0)\}) \geq 0.
\end{equation}
Choosing, in particular, $\alpha = \I$ this implies
$$\frac{1}{2\pi}\:\frac{2q_1(\I)}{2\I} + \frac{1}{4\pi^2}\nu(\{(0,0)\}) \geq 0.$$
But this is only possible if $q_1(\I) = d_1 \I$ for some $d_1 \geq 0$. The left hand side of  \myref{eq563} now takes the form
$$\frac{1}{2\pi}\:\frac{q_1(\alpha)+ d_1\I}{2\I} - \frac{\I}{4\pi^2}\frac{1}{2}(\alpha+\I)\nu(\{(0,0)\}).$$
This is a holomorphic function in the variable $\alpha \in \C^+$ which is real valued and hence constant. This implies that
$$q_1(\alpha) = 2\pi\I b_1 - \frac{1}{2\pi}(\alpha + \I)\nu(\{(0,0)\})$$
for some $b_1 \in \R$. Recall that by definition  $q_1$ is a Herglotz-Nevanlinna function without a linear term, i.e.  the number $b$ from representation \myref{eqthm21} is equal to 0, hence by Proposition \ref{prop22} it holds that 
$$\lim\limits_{\alpha \hat{\to} \infty}\frac{q_1(\alpha)}{\alpha} = 0.$$
 This implies that $\nu(\{(0,0)\}) = 0$ and $q_1(\alpha) = 2\pi\I\:b_1$ where $b_1 = \tfrac{d_1}{2\pi} \geq 0$.

In the same way we fix now the first variable  $\beta \in \C^+$, consider  the function $\til{q}_2$ defined as $\til{q}_2(z) := q(\beta,z)$ for $z \in \C^+$. The same reasoning  gives that $q_2(\beta) = 2\pi\I\:b_2$ for some $b_2 \geq 0$. 

Returning to representation \myref{eq52} we see that the second term is equal to 0, the third term becomes $b_2 z_2$ while the fourth term becomes $b_1z_1$. This then completes representation \myref{eqthm51}. 

It remains to show that the measure $\mu$ satisfies condition \myref{eqthm53}. We begin by recalling that $\nu$ satisfies condition \myref{04eq5} which implies that  also 
$$\iint_{[0,2\pi)^2}\sum_{(n_1,n_2) \in \N^2}\bar{w_1}^{n_1}w_2^{n_2}\E^{\I n_1 s_1}\E^{-\I n_2 s_2}\diff\nu(s_1,s_1) \equiv 0$$
for any $(w_1,w_2) \in \D^2$ since a geometric series permits the interchange of integration and summation. An analogous statement holds also for the conjugate of the above series. We thus conclude that
\begin{equation}\label{eq572}
\iint_{[0,2\pi)^2}\Re\left[\frac{\bar{w_1}w_2\E^{\I s_1}\E^{-\I s_2}}{(1-\bar{w_1}\E^{\I s_1})(1-w_2\E^{-\I s_2})}\right]\diff\nu(s_1,s_2) \equiv 0.
\end{equation}

The proof of Theorem \ref{thm5} so far will allow us to change the area of integration in  \myref{eq572} into the open square. To this end, we  begin by splitting the area of integration in formula \myref{eq572} into four parts as previously. The integral over the set $\{0\} \times \{0\}$ vanishes since we have shown that $\nu(\{(0,0)\}) = 0$. 

The integrals over the sets $\{0\} \times (0,2\pi)$ and $(0,2\pi) \times \{0\}$ are also equal to $0$. To see this  recall the functions $q_1$ and $q_2$, which were Herglotz-Nevanlinna functions in one variable defined via integrals over these lines, have been shown to be identically equal to $d_1 \I$ and $d_2 \I$ respectively, where $d_1,d_2 \geq 0$. Since the measure appearing in the representation of a Herglotz-Nevanlinna function in one variable is unique the measures $\mu_1$ and $\mu_2$ have to be equal to $d_1 \lambda_\R$ and $d_2 \lambda_\R$ respectively where $\lambda_\R$ denotes the Lebesgue measure on $\R$. Integration over the sets $\{0\} \times (0,2\pi)$ and $(0,2\pi) \times \{0\}$ in formula \myref{eq572} thus reduces to integrals of the from
$$\int_{[0,2\pi)}\E^{\I n s}\diff s$$
which vanish for  $n \in \Z\setminus\{0\}$.

Finally, the integral over the open square $(0,2\pi)^2$ remains and gives
\begin{equation}\label{eq573}
\iint_{(0,2\pi)^2}\Re\left[\frac{\bar{w_1}w_2\E^{\I s_1}\E^{-\I s_2}}{(1-\bar{w_1}\E^{\I s_1})(1-w_2\E^{-\I s_2})}\right]\diff\nu(s_1,s_2) \equiv 0.
\end{equation}

We can now change the area of integration in formula \myref{eq573} to $\R^2$ with the same change of coordinates used throughout the proof and expressed with the function $\varphi$ from Figure \ref{fig1}.  The integrand then transforms as
$$\frac{\bar{w_1}w_2\E^{\I s_1}\E^{-\I s_2}}{(1-\bar{w_1}\E^{\I s_1})(1-w_2\E^{-\I s_2})} = \frac{(\bar{z_1} + \I)(z_2-\I)(t_1-\I)(t_2+\I)}{4(t_1 - \bar{z_1})(t_2-z_2)}.$$
We also get a factor
$$\frac{4}{(1+t_1^2)(1+t_2^2)}$$
that comes from $\diff\varphi$. Formula \myref{eq573} thus transforms into
$$\iint_{\R^2}\Re\left[\frac{(\bar{z_1} + \I)(z_2-\I)(t_1-\I)(t_2+\I)}{4(t_1 - \bar{z_1})(t_2-z_2)}\right]\frac{4}{(1+t_1^2)(1+t_2^2)}\diff\mu(t_1,t_2) \equiv 0.$$
Since
\begin{multline*}
\frac{(\bar{z_1} + \I)(z_2-\I)(t_1-\I)(t_2+\I)}{4(t_1 - \bar{z_1})(t_2-z_2)}\cdot\frac{4}{(1+t_1^2)(1+t_2^2)}  \\
= \left(\frac{1}{t_1-\bar{z_1}}-\frac{1}{t_1+\I}\right)\left(\frac{1}{t_2-z_2}-\frac{1}{t_2-\I}\right),
\end{multline*}
this implies that the measure $\mu$ does indeed satisfy condition \myref{eqthm53}. We have thus proven that every Herglotz-Nevanlinna function in two variables admits a representation of the form \myref{eqthm51}. 

Conversely, let $q$ be a function defined on $\C^{+2}$ by \myref{eqthm51} with the number $a,b_1,b_2$ and the measure $\mu$ satisfying all the listed properties. The integral that appears in representation \myref{eqthm51} is a well-defined expression since the measure $\mu$ satisfies condition \myref{eqthm52}. It is then easy to see that a function $q$ defined in this way is holomorphic on $\C^{+2}$ since the kernel $K$ is holomorphic and locally uniformly bounded on compact subsets of $\C^{+2}$. 

To see that also $\Im[q] \geq 0$ consider the imaginary part of $q$ given by  \myref{eqthm51}, which is
$$\Im[q(z_1,z_2)] = b_1\Im[z_1] + b_2\Im[z_2] + \frac{1}{\pi^2}\iint_{\R^2}\Im[K\big((z_1,z_2),(t_1,t_2)\big)]\diff\mu(t_1,t_2).$$
Note that we are allowed to move the imaginary part into the integral due to $\mu$ being a real measure. It is now obvious that the first two terms are non-negative. 

To see that the third term is also non-negative observe that
\begin{multline*}
\Im[K\big((z_1,z_2),(t_1,t_2)\big)] = \pois\big((z_1,z_2),(t_1,t_2)\big) \\
-\frac{1}{2}\Re\left[\left(\frac{1}{t_1-\bar{z_1}}-\frac{1}{t_1+\I}\right)\left(\frac{1}{t_2-z_2}-\frac{1}{t_2-\I}\right)\right].
\end{multline*}
Since the measure $\mu$ satisfies property \myref{eq573} we have that
$$\iint_{\R^2}\Im[K]\,\diff\mu = \iint_{\R^2}\pois\,\diff\mu \geq 0$$
where the last inequality comes from the positivity of the Poission kernel and the dependence on the variables has been suppressed to shorten notation. This finishes the proof.
\endproof

\begin{opomba1}\label{vlad}
As mentioned in the introduction an integral representation of the form \eqref{eqthm51} already appears in \cite{vladimirov1} even for the case of $n \geq 2$ variables. However, it is only shown - by a completely different method - that every Herglotz-Nevanlinna function admits such a representation, but it is a priori assumed that the measure is the boundary measure of the representing function. It is  not discussed which measures actually can appear there. 

In \cite{vladimirov2} the authors use also a similar change of variables as in the present paper in order to find a characterization of Herglotz-Nevanlinna functions. However, this representation is not as simple as  \eqref{eqthm51}. Basically, all the integrals that come from the boundary of the area of integration are still present and hence also the corresponding Nevanlinna-condition is much more involved. 
\end{opomba1}

For convenience we highlight some minor results that appeared within in the proof of Theorem \ref{thm5}.

\begin{posledica1}\label{posfour}
Let $q$ be a Herglotz-Nevanlinna function in two variables. Then the following four statements hold.
\begin{itemize}
\item[(i)]{If there exists a point $(\zeta_1,\zeta_2) \in \C^{+2}$ such that $\Im[q(\zeta_1,\zeta_2)] = 0$ then $q(z_1,z_2) \equiv q(\zeta_1,\zeta_2)$ for all $(z_1,z_2) \in \C^{+2}$.}
\item[(ii)]{The imaginary part of $q$ can be represented as
$$\Im[q(z_1,z_2)] = b_1\Im[z_1] + b_2\Im[z_2] + \frac{1}{\pi^2}\iint_{\R^2}\pois\big((z_1,z_2),(t_1,t_2)\big)\diff\mu(t_1,t_2)$$
where $b_1,b_2$ and $\mu$ are as in Theorem \ref{thm5}.}
\item[(iii)]{The number $a$ from Theorem \ref{thm5} is equal to
$$a = \Re[q(\I,\I)].$$
}
\item[(iv)]{For every $\alpha,\,\beta \in \C^{+}$ it holds that
$$b_1 = \lim\limits_{z \hat{\to} \infty}\frac{q(z,\alpha)}{z} \quad \och \quad b_2 = \lim\limits_{z \hat{\to} \infty}\frac{q(\beta,z)}{z}$$
where $b_1,b_2$ are as in Theorem \ref{thm5}, in particular, the limits are  independent of $\alpha$ and $\beta$, respectively.}
\end{itemize}
\end{posledica1}

A further implication of Corollary \ref{posfour} is given by the following statement.

\begin{posledica1}\label{pos7}
Let $q$ be a Herglotz-Nevanlinna function in two variables and let $\alpha,\,\beta \in \C^{+}$ be arbitrary. Then there exist constants $c_1,c_2 \leq 0$ independent of $\alpha$ and $\beta$ such that
$$c_1 = \lim\limits_{z \hat{\to} 0}z\:q(z,\alpha) \quad \och \quad c_2 = \lim\limits_{z \hat{\to} 0}z\:q(\beta,z).$$
\end{posledica1}

\proof
Applying the change of variables $z \mapsto -\tfrac{1}{z}$ leads to Herglotz-Nevanlinna functions $z \mapsto q(-\tfrac{1}{z},\alpha)$ and $z \mapsto q(\beta,-\tfrac{1}{z})$ for which  Corollary \ref{posfour}(iv)  implies the claim.
\endproof

The proof of Theorem \ref{thm5} has also given us additional information about measures satisfying condition \myref{04eq5}.

\begin{posledica1}\label{pos5}
Let $\nu$ be a finite positive Borel measure on $[0,2\pi)^2$ satisfying condition \myref{04eq5}. Then $\nu(\{(0,0)\}) = 0$ and there exist constants $e_1,e_2 \geq 0$ such that $\nu|_{\{0\} \times (0,2\pi)} = e_1\lambda_{(0,2\pi)}$ and  $\nu|_{(0,2\pi) \times \{0\}} = e_2\lambda_{(0,2\pi)}$. In particular it holds that
$$\iint_{(0,2\pi)^2}\E^{\I m_1 s_1}\E^{\I m_2 s_2}\diff\nu(s_1,s_2) = 0$$
for every pair of indices $m_1,m_2 \in \Z$ satisfying $m_1m_2 < 0$.
\end{posledica1}

We finish this section with some examples of representations of Herglotz-Nevanlinna functions in two variables.

\begin{primer1}\label{ex1}
Let
$$q(z_1,z_2) = - \frac{1}{z_2} .$$
Then $q$ is a Herglotz-Nevanlinna function in two variables which can easily be shown by a direct computation of its imaginary part. Corollary \ref{posfour} now says that
$$a = \Re[q(\I,\I)] = 0$$
while choosing $\alpha=\beta  = \I$ we get that
$$b_1 = \lim\limits_{z \hat{\to} \infty}\frac{q(z,\I)}{z} = 0 \quad \och \quad b_2 = \lim\limits_{z \hat{\to} \infty}\frac{q(\I,z)}{z} = 0.$$
The measure $\mu$ can also be reconstructed using proposition \ref{pos4} and is equal to
$$\mu = \lambda_\R \otimes \pi\delta_0.$$
Note that if $q$ is regarded as a function in just one variable the representing measure (in Theorem \ref{thm21}) is only the Dirac measure $\delta_0.$
\end{primer1}

\begin{primer1}
Let
$$q(z_1,z_2) = 2 + z_1 + \frac{z_1z_2 + z_2 - z_1 - 1}{z_1 + z_2}.$$
Then $q$ is a Herglotz-Nevanlinna function in two variables with
$$a = \Re[q(\I,\I)] = 2,\quad  b_1 = \lim\limits_{z \hat{\to} \infty}\frac{q(z,\I)}{z} = 1, \quad  b_2 = \lim\limits_{z \hat{\to} \infty}\frac{q(\I,z)}{z} = 0,$$
and the measure $\mu$   equals 
$$\mu = \pi g \chi_{\{t_1 = -t_2\}}\lambda_{\R^2}$$
where the function $g$ is defined as $g(t_1,t_1) = -t_1t_2 - t_2 + t_1 + 1$.
\end{primer1}

\begin{primer1}
Let
$$q(z_1,z_2) = 1 + (2 + \sqrt{z_1})(3 + \sqrt{z_2})$$
where the branch cut of the square root function is taken along the negative real line. Then $q$ is a Herglotz-Nevanlinna function in two variables which can again be shown by  a direct computation of its imaginary part. We then have
$$a = \Re[q(\I,\I)] = 7+\frac{5}{\sqrt{2}}, \quad b_1 = \lim\limits_{z \hat{\to} \infty}\frac{q(z,\I)}{z} = 0,   \quad b_2 = \lim\limits_{z \hat{\to} \infty}\frac{q(\I,z)}{z} = 0,$$
and the measure $\mu$ is equal to
$$\mu = 3 \pi h_1 \chi_{\{t_1 < 0\}}\lambda_\R \otimes \lambda_\R + \lambda_\R \otimes 2 \pi h_1\chi_{\{t_2 < 0\}}\lambda_\R + h_2\chi_{\{t_1t_2 < 0\}}\lambda_{\R^2}$$
where the function $h_1$ is defined as $h_1(t) = \sqrt{-t}$ and the function $h_2$ is defined as $h_2(t_1,t_2) = \sqrt{-t_1t_2}$.
\end{primer1}

\section{Properties of the representing measure}\label{measure4}

We now return to the question of describing the measure $\mu$ in terms of the function $q$. In the one-variable case this is done via the classic Stieltjes inversion formula \cite{kac}. Here we present an elementary two-dimensional analog to this formula.

\begin{trditev1}\label{pos4}
Let $q$ be a Herglotz-Nevanlinna function in two variables. Then the measure $\mu$  in representation \myref{eqthm51}  is unique and can be determined from the boundary values of $q$. 
More precisely, let $\psi\colon \R^2 \to \R^2$ be a $\cont^1$ function such that
$$|\psi(x_1,x_2)| \leq \frac{C}{(1+x_1^2)(1+x_2^2)}$$
for some constant $C \geq 0$ and all $(x_1,x_2) \in \R^2$. Then
\begin{equation}\label{still}
\iint_{\R^2}\psi(t_1,t_2)\diff\mu(t_1,t_2)=\lim\limits_{\substack{y_1\:\downarrow\:0 \\ y_2\:\downarrow\:0}}\iint_{\R^2}\psi(x_1,x_2)\Im[q(x_1 + \I\:y_1,x_2+\I\:y_2)]\diff x_1\diff x_2 .
\end{equation}
\end{trditev1}

\proof
We begin by using statement (ii) of Corollary \ref{posfour} to rewrite the right-hand side of equality \myref{still} as
\begin{multline*}
\lim\limits_{\substack{y_1\:\downarrow\:0 \\ y_2\:\downarrow\:0}}\left.\iint_{\R^2}\psi(x_1,x_2)\right(b_1y_1+b_2y_2 \\
+ \left.\frac{1}{\pi^2}\iint_{\R^2}\pois\big((x_1+ \I\:y_1,x_2 + \I\:y_2),(t_1,t_2)\big)\diff\mu(t_1,t_2)\right)\diff x_1\diff x_2.
\end{multline*}
The part involving the term $b_1y_1+b_2y_2$ is equal to $0$ since Lebesgue's dominated convergence theorem allows us to change the order of the limit and the integral. What remains is the part involving the Poisson kernel where we can use Fubini's theorem to change the order of integration. Another application of Lebesgue's dominated convergence theorem allows us to change the order of the limit and the first integral. We thus arrive at
$$\frac{1}{\pi^2}\iint_{\R^2}\lim\limits_{\substack{y_1\:\downarrow\:0 \\ y_2\:\downarrow\:0}}\iint_{\R^2}\psi(x_1,x_2)\pois\big((x_1+ \I\:y_1,x_2 + \I\:y_2),(t_1,t_2)\big)\diff x_1\diff x_2\diff\mu(t_1,t_2).$$
It remains to observe that by a well known property of the Poisson kernel the inner integral equals 
 $$\lim\limits_{\substack{y_1\:\downarrow\:0 \\ y_2\:\downarrow\:0}}\iint_{\R^2}\psi(x_1,x_2)\pois\big((x_1+ \I\:y_1,x_2 + \I\:y_2),(t_1,t_2)\big)\diff x_1\diff x_2 = \pi^2\psi(t_1,t_2).$$
 
In order to show the uniqueness of the representing measure, suppose that  representation \myref{eqthm51} for the function $q$ holds for some measures $\mu_1$ and $\mu_2$. Recall that Corollary \ref{posfour} shows that the numbers $a,b_1,b_2$ are uniquely determined by $q$. Using  \myref{still} we see that the left-hand side in this formula   is the same for both $\mu_1$ and $\mu_2$. This implies
$$\iint_{\R^2}\psi(t_1,t_2)\diff\mu_1(t_1,t_2) = \iint_{\R^2}\psi(t_1,t_2)\diff\mu_2(t_1,t_2)$$
for all functions $\psi$ as above, which   is   possible only  if $\mu_1 \equiv \mu_2$.
\endproof

Recall that in the second part of the proof of Theorem \ref{thm5} we only required the growth condition to show that the integral involving the kernel function $K$ is well defined while the Nevanlinna condition is needed only show that the integral of $\Im[K]$ is non-negative. We illustrate this by the following example of a finite measure. 

\begin{primer1}
Considering the function defined by representation \myref{eqthm51} with $a = 0$, $b_1 = 0$, $b_2 = 0$ and $\mu = \pi^2\delta_{(0,0)}$, it is given by 
$$\iint_{\R^2}K\big((z_1,z_2),(t_1,t_2)\big)\diff\delta_{(0,0)}(t_1,t_2) = \frac{\I(z_1+\I)(z_2+\I)}{2z_1z_2}-\I.$$
Note that the measure $\pi^2\delta_{(0,0)}$ does not satisfy the Nevanlinna condition \myref{eqthm53} and hence the above function is not a Herglotz-Nevanlinna function.
\end{primer1}

We show now that finite measures actually cannot satisfy the Nevanlinna condition.
\begin{trditev1}\label{pos9}
Let $q$ be a Herglotz-Nevanlinna function in two variables and let $\mu$ be the representing measure. Then $\mu$ cannot be a finite measure unless it is identically equal to $0$.
\end{trditev1}

\proof
Recall first that for a Herglotz-Nevanlinna function $\til{q}$ of one variable with $b = 0$  and representing measure $\til{\mu}$  it holds that
$$\lim\limits_{y \to \infty}y\:\Im[\til{q}(\I\:y)] = \frac{1}{\pi}\int_{-\infty}^{\infty}\diff\til{\mu}(t).$$
 Observe that this identity holds even if one (and thus both) sides are equal to $+\infty$.

For a Herglotz-Nevanlinna function $q$ in two variables that has $b_1 = b_2 = 0$, where the numbers $b_1,b_2$ are as in Theorem \ref{thm5},  Corollary \ref{posfour}(ii)  and Lebesgue's monotone convergence theorem imply  that
\begin{equation}\label{finite}
\lim\limits_{y \to \infty}y^2\:\Im[q(\I\:y,\I\:y)] = \frac{1}{\pi^2}\iint_{\R^2}\diff\mu(t_1,t_2).
\end{equation}
 As in the one-variable case the identity remains valid if one (and thus both) sides are equal to $+\infty$.

Suppose now that $ q$ is a Herglotz-Nevanlinna function in two variables with a finite representing measure $\mu$. Note that  $\mu$ is also a representing measure for the Herglotz-Nevanlinna function $\hat{q}$ defined as $\hat{q}(z_1,z_2): = q(z_1,z_2) - a -b_1z_1 - b_2z_2$, where  $a,b_1,b_2$ are  as in Theorem \ref{thm5}.

The finiteness of  $\mu$ and equation \eqref{finite} imply $\lim_{y \to \infty}y^2\:\Im[\hat{q}(\I\:y,\I\:y)] < \infty$ and hence  $\lim_{y \to \infty}y\:\Im[\hat{q}(\I\:y,\I\:y)] = 0$. In terms of the Herglotz-Nevanlinna function $\til{q}$ defined as $\til{q}(z) := \hat{q}(z,z)$ this translates to
$$\frac{1}{\pi}\int_{-\infty}^{\infty}\diff\til{\mu}(t) = 0.$$
In follows that $\til{q} \equiv 0$ and so $\hat{q} \equiv 0$ on the diagonal in $\C^{+2}$. In particular $\hat{q}(\I,\I) = 0$ which together with statement (i) of Corollary \ref{posfour} implies $\hat{q} \equiv 0$ and thus $\Im[\hat{q}] \equiv 0$. Statement (ii) of Corollary \ref{posfour} now show that $\mu \equiv 0$ is the only possibility.
\endproof

We can now in fact say even more about the measures that are allowed in representation \myref{eqthm51}, namely that   points in $ \R^2$ are always zero sets.

\begin{trditev1}\label{pos10}
Let $q$ be a Herglotz-Nevanlinna function in two variables and let $\mu$ be the representing. Then $\mu(\{(t_{01},t_{02})\}) = 0$ for any point $(t_{01},t_{02}) \in \R^2$.
\end{trditev1}

\proof
We begin by observing that
$$\lim\limits_{\substack{z_1\:\hat{\to}\:t_{01} \\ z_2\:\hat{\to}\:t_{02}}}(z_1-t_{01})(z_2-t_{02})K\big((z_1,z_2),(t_1,t_2)\big) = -\frac{\I}{2}\chi_{\{(t_{01},t_{02})\}}(t_1,t_2)$$
for any point $(t_{01},t_{02}) \in \R^2$.

For any fixed $\beta \in \C^+$ we calculate also that
$$\lim\limits_{z_2 \hat{\to} t_{02}}(z_2-t_{02})q(\beta,z_2)  =  -\lim\limits_{\omega_2 \hat{\to} \infty}\tfrac{1}{\omega_2}\til{q_2}(\beta,\omega_2) = - \til{b_2}(t_{02}).$$
Here we used the variable change $z_2-t_{02} = -\tfrac{1}{\omega_2}$ along with statement (iv) of Corollary \ref{posfour} for the Herglotz-Nevanlinna function 
$$\til{q}_2\colon (\omega_1,\omega_2) \mapsto q(\omega_1,-\tfrac{1}{\omega_2}+t_{02}).$$
Note that the number $\til{b_2}(t_{02})$ does of course depend on $t_{02}$ but it does not depend on $\beta \in \C^+$. This implies that
$$\lim\limits_{\substack{z_1\:\hat{\to}\:t_{01} \\ z_2\:\hat{\to}\:t_{02}}}(z_1-t_{01})(z_2-t_{02})q(z_1,z_2) = -\lim\limits_{z_1\:\hat{\to}\:t_{01}}(z_1-t_{01})\til{b_2}(t_{02}) = 0.$$

On the other hand we can use Theorem \ref{thm5}, Lebesgue's dominated convergence theorem and our starting observation to show that
$$\lim\limits_{\substack{z_1\:\hat{\to}\:t_{01} \\ z_2\:\hat{\to}\:t_{02}}}(z_1-t_{01})(z_2-t_{02})q(z_1,z_2) = -\frac{\I}{2\pi^ 2}\mu(\{(t_{01},t_{02})\}).$$
This finishes the proof.
\endproof

\section*{Conclusion}

Theorem \ref{thm5} provides the anticipated generalization of Theorem \ref{thm21} to the case of two variables while improving the previous results of Vladimirov \cite{vladimirov1,vladimirov2} as discussed in the introduction and in Remark \ref{vlad}. A similar type of improvement of Valdimirov's results for the general case $n > 2$ will be considered in an upcoming work.

The theorem has also allowed us to understand that the divide between the cases $n=1$ and $n=2$ exists first and foremost in the class of measure that can represent a Herglotz-Nevanlinna function in the respected dimension. While the one dimensional case has been completely understood since the appearance of Theorem \ref{thm21} about a century ago we have now seen for example that all representing measure of Herglotz-Nevanlinna functions in two variables are atomless. The properties of the corresponding class of measures for the case $n > 2$ will be considered along the upcoming generalization of Theorem \ref{thm5}.

\section*{Acknowledgement}

The authors thank Ragnar Sigurdsson for interesting discussions on the subject and careful reading of the manuscript. 

\bibliographystyle{amsplain}

\begin{thebibliography}{99}

\bibitem{agler2}
J. Agler, J. E. McCarthy, N. J. Young, \emph{Operator monotone functions and L\"{o}wner functions of several variables}, Ann. of Math. (2) \textbf{176} (2012), no. 3, 1783--1826. 

\bibitem{agler}
J. Agler, R. Tully-Doyle, N. J. Young, \emph{Nevanlinna representations in several variables}, arXiv:1203.2261v2.

\bibitem{cauer}
W. Cauer, \emph{The Poisson integral for functions with positive real part}, Bull. Amer. Math. Soc. \textbf{38} (1932), no. 10, 713--717.

\bibitem{vladimirov3}
Yu. N. Drozhzhinov and V. S. Vladimirov, \emph{Holomorphic functions in a polydisc with non-negative imaginary part}, Mat. Zametki \textbf{15} (1974), 55--61. English transl. in Math. Notes \textbf{15} (1974), 31--34.

\bibitem{kac}
I. S. Kac and M. G. Krein, \emph{R-functions--analytic functions mapping the upper half-plane into itself}, Amer. Math. Soc. Transl. (2) \textbf{103} (1974), 1--18.

\bibitem{koranyi}
A.~Kor\'{a}nyi and L. Puk\'{a}nszky, \emph{Holomorphic function with positive real part on polycylinders}, Trans. Amer. Math. Soc. \textbf{108} (1963), 449--456.

\bibitem{nevan1}
R. Nevanlinna, \emph{Asymptotische Entwicklungen beschr\"{a}nkter Funktionen und das Stieltjessche Momentenproblem}, Ann. Acad. Sci. Fenn. (A) \textbf{18 (5)} (1922), 1--53.

\bibitem{vladimirov1}
V. S. Vladimirov, \emph{Generalized function in mathematical physics},  "Nauka", 1979; English transl., Mir publishers, Moscow, 1979.
 
\bibitem{vladimirov2}
V. S. Vladimirov, \emph{Holomorphic functions with non-negative imaginary part in a tubular domain over a cone},  Mat. Sb. \textbf{79 (121)} (1969), 182--152. English transl. in Math. USSR-Sb. \textbf{8} (1969), 125--146.

\end{thebibliography}

\end{document}